\documentclass[JEDP]{cedram-sem}
\usepackage{amsfonts}
\usepackage{amsthm}
\usepackage{a4wide}
\usepackage[utf8]{inputenc}
\usepackage{amsmath}
  \usepackage{xcolor} 
\usepackage{amssymb}
 \usepackage{subfigure}
\usepackage[mathscr]{eucal}
\usepackage{dsfont}
\usepackage{pstricks}
\usepackage{graphicx}
\usepackage[colorlinks=true, linkcolor=purple, citecolor=purple]{hyperref}

\newtheorem{theorem}{Theorem}[section]

\newtheorem{proposition}[theorem]{Proposition}
\newtheorem{corollary}[theorem]{Corollary}

\newtheorem{ques}{Question}

\newcommand{\nwc}{\newcommand}
\nwc{\eps}{\epsilon}
\nwc{\ep}{\epsilon}
\nwc{\vareps}{\varepsilon}
\nwc{\Oph}{\operatorname{Op}_\hbar}
\nwc{\la}{\langle}
\nwc{\ra}{\rangle}

\nwc{\mf}{\mathbf} %Latex (as in \bf not tilted math letters)
\nwc{\blds}{\boldsymbol} %Latex
\nwc{\ml}{\mathcal} %Latex

\nwc{\defeq}{\stackrel{\rm{def}}{=}}

\nwc{\cE}{\ml{E}}
\nwc{\cN}{\ml{N}}
\nwc{\cO}{\ml{O}}
\nwc{\cP}{\ml{P}}
\nwc{\cU}{\ml{U}}
\nwc{\cV}{\ml{V}}
\nwc{\cW}{\ml{W}}
\nwc{\M}{\ml{M}}
\let \H \relax
\nwc{\H}{\ml{H}}

\nwc{\tU}{\widetilde{U}}
\nwc{\IN}{\mathbb{N}}
\nwc{\IR}{\mathbb{R}}
\nwc{\R}{\mathbb{R}}
\nwc{\IZ}{\mathbb{Z}}
\nwc{\Z}{\mathbb{Z}}
\nwc{\N}{\mathbb{N}}
\nwc{\IC}{\mathbb{C}}
\nwc{\C}{\mathbb{C}}
\nwc{\IT}{\mathbb{T}}
\nwc{\T}{\mathbb{T}}
\nwc{\IS}{\mathbb{S}}
\nwc{\tP}{\widetilde{P}}
\nwc{\tPi}{\widetilde{\Pi}}
\nwc{\tV}{\widetilde{V}}
\nwc{\supp}{\operatorname{supp}}
\nwc{\rest}{\restriction}
\let \d \relax
\nwc{\d}{\partial}
\nwc{\Cor}{\mathscr{C}}
\nwc{\F}{\mathscr{F}}

\nwc{\CLam}{\overline{\C}_+^\Lambda}

\nwc{\todo}[1]{$\clubsuit$ {\tt #1}}

\DeclareMathOperator{\Vol}{Vol}

\DeclareMathOperator{\Res}{Res}
\DeclareMathOperator{\Sp}{Sp}

\DeclareMathOperator{\Conv}{Conv}
\DeclareMathOperator{\singsupp}{Sing\,supp}

\renewcommand{\Re}{\operatorname{Re}}
\let \div \relax \newcommand{\div}{\operatorname{div}}

\title 
[Length orthospectrum on flat tori]
{Length orthospectrum and the correlation function on flat tori}

\author
[N. V. \lastname{Dang}]
{\firstname{Nguyen Viet} \lastname{Dang}}
\address{Sorbonne Universit\'e, IMJ-PRG, 75252 Paris Cedex 05, France}
\email{dang@imj-prg.fr}

\author
[M. \lastname{L\'eautaud}]
{\firstname{Matthieu} \lastname{L\'eautaud}}
\address{Laboratoire de Math\'ematiques d'Orsay, UMR 8628, Universit\'e Paris-Saclay, CNRS, B\^atiment 307, 91405 Orsay Cedex France}
\email{matthieu.leautaud@universite-paris-saclay.fr}

\author
[G. \lastname{Rivi\`ere}]
{\firstname{Gabriel} \lastname{Rivi\`ere}}
\address{Laboratoire de Math\'ematiques Jean Leray, Nantes Universit\'e, UMR CNRS 6629, 2 rue de la Houssini\`ere, 44322 Nantes Cedex 03, France}
\address{Institut Universitaire de France, Paris, France}%\thanks{J.~D. is supported by .}
\email{gabriel.riviere@univ-nantes.fr}

\begin{document}

\begin{abstract}
This note presents some of the results obtained in~\cite{DangLeautaudRiviere21} and it has been
the object of a talk of the second author during the Journées “Équations aux
Dérivées Partielles” (Obernai, june 2022).
We study properties of geodesics that are orthogonal to two convex subsets of the flat torus $\T^d$. We discuss meromorphic properties of a geometric Epstein zeta function associated to the set of lengths of such orthogeodesics.
We also define the associated length distribution and discuss singularities of its Fourier transform.
Our analysis relies on a fine study of the dynamical correlation function of the geodesic flow on the torus and the definition of anisotropic Sobolev spaces that are well-adapted to this integrable dynamics.
\end{abstract}

% French abstract
%\begin{altabstract}
%
%\end{altabstract}

\maketitle
%\tableofcontents

\setcounter{tocdepth}{2}
\tableofcontents

\section{Introduction: three questions on the torus}

In this note, we review some of the recent results obtained by the authors in~\cite{DangLeautaudRiviere21}, in the light of three different but related questions. The first topic concerns a family of geometric zeta functions associated to the set of lengths of orthogeodesics to two convex bodies on the torus. The second topic concerns Poisson formul{\ae} associated to the same set. And the last topic is related to decay of the dynamical correlation function for the geodesic flow on the torus.

%%%%%%%%%%%%%%%%%%%%%%%%%%%%%%
\subsection{Question 1: zeta functions}
\label{s:zeta-fct}
The original and paradigmatic zeta function is of course the Riemann zeta function: 
  $$\zeta_{\text{Riem}}(s):=\sum_{n \in \N^*}\frac{1}{n^{s}} .$$
Absolutely convergent and holomorphic on the half-plane $\Re(s) >1$, it is well-known that $\zeta_{\text{Riem}}$ continues meromorphically to the whole complex plane $\C$, with a single pole at $s=1$ whose residue is given by $\Res_{s=1}(\zeta_{\text{Riem}}) =1$.  The localization of its nontrivial zeroes is related to the growth of prime numbers and constitutes the famous Riemann hypothesis.
  
A closely related zeta function was introduced by Epstein~\cite{Epstein03}: given $q \in \R^d$, set
 $$\zeta_{\text{Eps}}(q,s):=\sum_{\xi\in\IZ^d\setminus\{-q\}}\frac{1}{|\xi+q|^s}, $$
 where $|\cdot|$ denotes the Euclidean norm in $\R^d$.
As for $\zeta_{\text{Riem}}$, this expression is absolutely convergent (and thus holomorphic) for $\Re(s)>d$. Epstein proved that $\zeta_{\text{Eps}}$ extends meromorphically to the whole complex plane $\C$, with a single pole at $s=d$ whose residue is given by $\Res_{s=d}(\zeta_{\text{Eps}}) =1$.

More generally, given a discrete set of real numbers (having at most polynomial growth), one can form from this set a zeta function and one hopes that meromorphic properties of the resulting zeta function will furnish useful information on the discrete set. Highly interesting discrete sets in geometry arise for instance from the lengths of certain geodesics on Riemannian manifolds. The latter allow to construct geometric zeta functions, for instance from the set of lengths of closed geodesics on negatively curved manifolds
(Selberg~\cite{Selberg56},
Smale~\cite{Smale67},
Ruelle~\cite{Ruelle76}, 
Fried~\cite{Fried95},
Rugh~\cite{Rugh96}, 
Giulietti Liverani Pollicott~\cite{GiuliettiLiveraniPollicott13}, 
Faure Tsujii~\cite{FaureTsujii13,FaureTsujii17}, Dyatlov Zworski~\cite{DyatlovZworski16}...)  or the set of  lengths of orthogeodesics on negatively curved manifolds, that is to say, geodesics that are orthogonal to two fixed (smooth) sets 
(Huber~\cite[Satz A]{Huber56},~\cite[Satz 2]{Huber59}, Margulis~\cite{Margulis69, Margulis04}, Parkkonen Paulin~\cite{ParkkonenPaulin2016},
Broise-Alamichel Parkkonen Paulin~\cite{BroiseParkkonenPaulin2019},
Dang Rivi\`ere~\cite{DangRiviere20d}...)
  
  With this geometric picture in mind, a first geometric zeta function on the torus $\T^d_1 := \R^d/ \Z^d$ (endowed with the flat metric) is the following two-point zeta function: given  $x_1,x_2 \in \T^d_1$, we set 
$$\mathcal{P}(x_1,x_2) = \{\text{nontrivial geodesic curves between } x_1 \text{ and }x_2 \}, \quad \ell(\gamma) = \text{length}(\gamma)>0 ,$$
 and define
 $$\zeta_{\text{Points}}(x_1,x_2,s):= \sum_{\gamma \in \mathcal{P}(x_1,x_2)}\frac{1}{\ell(\gamma)^s} , \quad \Re(s) \text{\ large enough} .$$
Parametrizing the set of such geodesics by $\xi \in \Z^d$, for instance by $\gamma_{\xi}(t) = (1-t)x_1 + t(x_2+\xi)$ for $t \in [0,1]$, we see that 
 $$\zeta_{\text{Points}}(x_1,x_2,s) 
 = 
 \sum_{\xi\in\IZ^d\setminus\{x_2-x_1\}}\frac{1}{|x_2-x_1 - \xi|^s} = \zeta_{\text{Eps}}(x_2-x_1,s) ,
$$
and this ``geometric zeta function'' coincides with the above described Epstein zeta function.

A richer geometric setting is the following. Let $K_1$ and $K_2$ be two strictly convex and compact subsets of $\IR^d$ ($d\geq 2$) with smooth boundaries $\partial K_1$ and $\partial K_2$. By strictly convex, we mean that the boundary of the convex set $K_i$ has all its sectional curvatures \emph{positive} (if $K$ is reduced to a point, then we adopt the convention that it is strictly convex with smooth boundary).  
 Through the canonical projection 
 \begin{align}
 \label{can-proj}
 \mathfrak{p}:\IR^d\rightarrow \IT^d=\IT_{2\pi}^d:=\IR^d/2\pi\IZ^d
 \end{align} (notice the slightly different normalization of the torus), the boundaries of $K_1$ and $K_2$ can be identified with immersed submanifolds of the flat torus (that may have selfintersection points). We fix an orientation on each submanifold $\partial K_i$ either by the outgoing normal vector to $K_i$ or by the ingoing one, which induces an orientation on $\Sigma_i:=\mathfrak{p}(\partial K_i)$. 
  Defining now $$\mathcal{P}_{K_1 ,K_2} =  \{\text{geodesic curves {\em directly} orthogonal to } K_1  \text{ and } K_2 \}, $$
the set of so-called orthogeodesics, we now set 
 $$\zeta_{\Conv}(K_1,K_2,s):= \sum_{\gamma \in \mathcal{P}_{K_1,K_2}}\frac{1}{\ell(\gamma)^s} , \quad \Re(s)  \text{\ large enough} .$$
In case $K_1,K_2$ are both points, we recover the definition of the function $\zeta_{\text{Points}}$.
 The same questions are in order for the function $s \mapsto \zeta_{\Conv}(K_1,K_2,s)$: 
 \begin{ques}
 \label{Q1}
 Where is $\zeta_{\Conv}(K_1,K_2,s)$ well-defined? Does it continue meromorphically to $\C$? Where are its poles and can we compute the associated residues?
\end{ques}
 
Question~\ref{Q1} is studied in Section~\ref{s:zeta} below.

\begin{figure}[ht]
\includegraphics[scale=0.3]{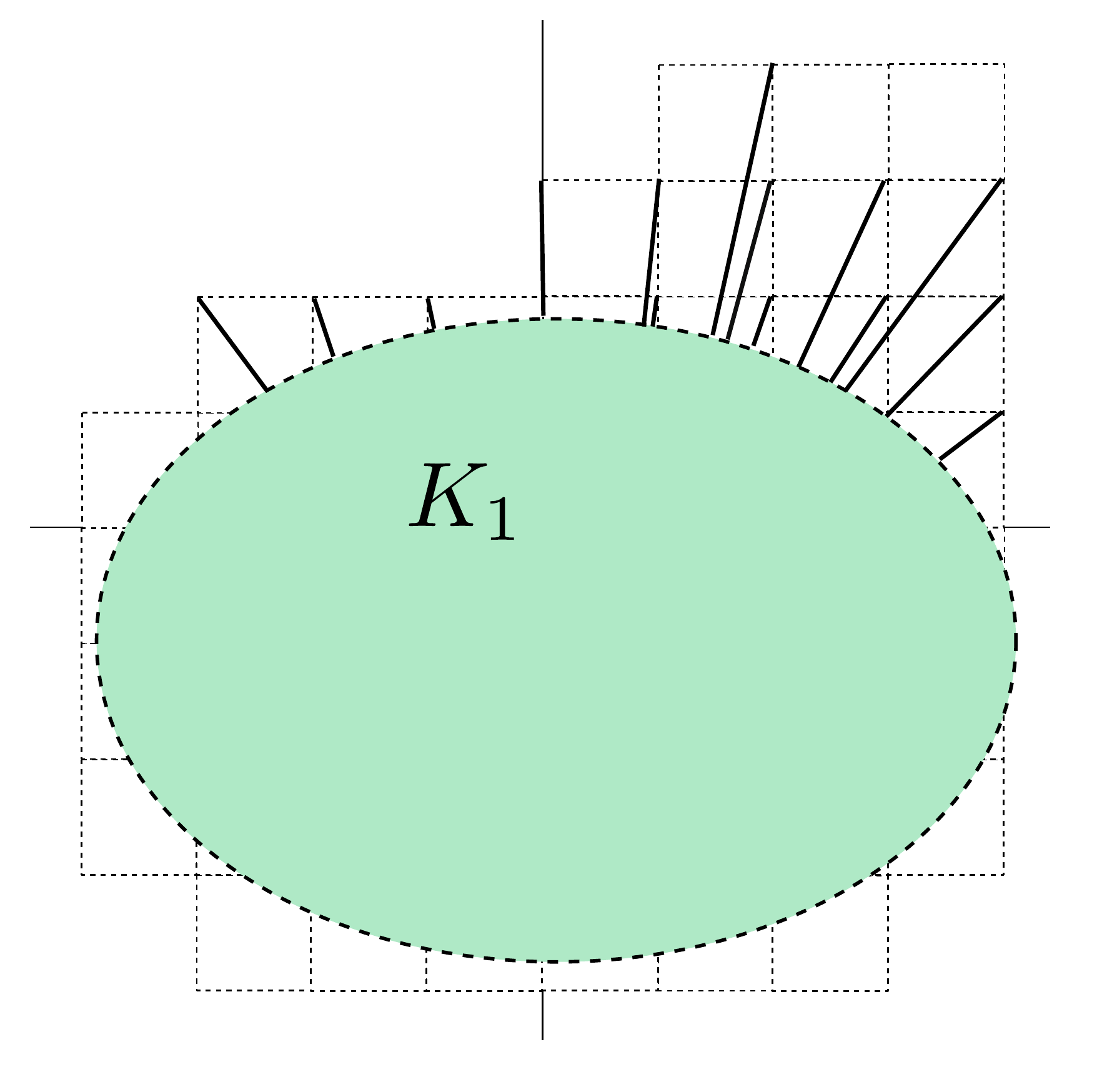}
\centering
\caption{\label{f:lattice}Lift of the orthogeodesic arcs when $K_2=\{0\}$.}
\end{figure}

\subsection{Question 2: Poisson formul\ae}
We start by recalling the famous Poisson formula, relating the values of any function $f$ on the lattice $\frac{2\pi}{\nu} \Z$, $\nu >0$ to those of its Fourier transform $\hat{f}$ on the dual lattice $\nu \Z$: for all $ f \in \mathcal{S}(\R)$
$$
\sum_{n \in \Z} \hat f (\nu n)= \frac{2\pi}{\nu} \sum_{k \in \Z} f\big(\frac{2\pi}{\nu} k\big) ,
$$
where $\hat f(\xi)  :=\F(f) (\xi)= \int_\R e^{-ix\xi}f(x)dx$. As usual, this formula may be rewritten as the Fourier transform of a Dirac comb:
\begin{equation}
\label{e:poisson}
\sum_{n \in \Z} e^{-i\nu n \tau}  = \F\left( \sum_{n \in \Z} \delta_{\nu n} \right) = \frac{2\pi}{\nu} \sum_{k \in \Z} \delta_{\frac{2\pi}{\nu} k} . 
\end{equation}
in the sense of temperate distributions, i.e. in  $\mathcal{S}'(\R)$. 
Now, given any discrete set of real numbers (growing at most polynomially), one can form a ``generalized Dirac comb'', and try to compute its Fourier transform (or, at least describe some of its properties). This would yield in principle a good candidate for a Poisson formula.

In particular,~\eqref{e:poisson} may be reinterpreted as follows: $\nu \Z$ might be read as $\nu \Z = \pm  \Sp(\sqrt{-\d_x^2})$ where $-\d_x^2$ is acting on functions on the torus $\T_{2\pi\nu}=\R/(2\pi\nu) \Z$, whereas $\frac{2\pi}{\nu} \Z$ might be read as the (symmetrized with respect to $0$) set of lengths of closed geodesics of the manifold $\T_{2\pi\nu}=\R/(2\pi\nu) \Z$.
A far-reaching generalization of the Poisson formul{\ae} was thus formulated in spectral theory by computing the Fourier transform of the Dirac comb 
$$ \mathsf{T}(t)=\sum_{\lambda\in \Sp(\sqrt{-\Delta_g})} \delta_{\lambda} \in \mathcal{S}^\prime(\mathbb{R}) ,$$ 
namely,
$$\widehat{\mathsf{T}}(\tau)=\sum_{\lambda\in \Sp(\sqrt{-\Delta_g})} e^{-i\tau \lambda }\in \mathcal{S}^\prime(\mathbb{R}) ,$$ 
where $\Delta_g$ is the Laplace-Beltrami operator on a compact Riemannian manifold and $\Sp(\sqrt{-\Delta_g})$ denotes the spectrum of its square root.
Convergence of the series in $\mathcal{S}'(\R)$ is due to the Weyl law, which ensures polynomial growth of the number of eigenvalues. The wave trace formula, proved by Chazarain~\cite{Chazarain74} and Duistermaat--Guillemin~\cite{DuistermaatGuillemin75} extending previous results by Selberg~\cite{Selberg56} and Colin de Verdi\`ere~\cite{ColindeVerdiere73}, states that the singular support of the distribution $\widehat{\mathsf{T}}$ is exactly the set of lengths of periodic geodesic curves for the metric $g$. Furthermore, when the geodesic flow is nondegenerate, it describes the singularity of $\widehat{\mathsf{T}}$ at each period in terms of geometric data attached to the periodic orbits and of distributions of the form $(t\pm \ell+i0)^{-1}$. In other words, the quantum spectrum determines the classical length spectrum and these wave trace formulas are often referred as generalized Poisson formulae. 
According to~\cite[p.72]{Hormander90}, the singularities in this formula can be rewritten as
\begin{equation}\label{e:finite-part}\lim_{y\rightarrow 0^+}\frac{1}{t\pm \ell +iy}=(t\pm \ell+i0)^{-1}=\text{FP}\left(\frac{1}{t\pm\ell}\right)-i\pi\delta_{\mp \ell},\end{equation}
where $\text{FP}\left(.\right)$ is the finite part of the (non-integrable) function $(t\pm \ell)^{-1}$.
 The contrast with~\eqref{e:poisson} comes from the fact the classical Dirac comb is symmetrized w.r.t. $0$, so that its Fourier transform only contains Dirac-type singularities.

\bigskip
Coming back to the geometric context of Section~\ref{s:zeta-fct}, we may define a geometric Dirac comb associated to the set of length of orthogeodesics to two strictly convex sets $K_1,K_2$ on $\T^d$ by
  \begin{align}
  \label{d:TKK}
  \mathsf{T}_{K_1,K_2}  = \sum_{\gamma \in \mathcal{P}_{K_1,K_2}} \delta_{\ell(\gamma)}   ,
  \end{align}
which converges in $\mathcal{S}'(\R)$ according to Proposition~\ref{l:starting-point}. The following questions are in order.
 \begin{ques}
 \label{Q2}
 What can be said about $$\F\left(\mathsf{T}_{K_1,K_2}\right)(\tau) =  \sum_{\gamma \in \mathcal{P}_{K_1,K_2}}e^{-i\tau \ell(\gamma) } ?$$ 
 Can one locate its singularities? Is $\F\left(\mathsf{T}_{K_1,K_2}\right)$ a Dirac comb?
 \end{ques}
 Question~\ref{Q2} is discussed in Section~\ref{s:poisson-ortho} below.

   \subsection{Question 3: the dynamical correlation function}
The last question we consider in these notes might seem mildly related with the previous ones at first sight. It arises both from hyperbolic dynamical systems and from the study of kinetic partial differential equations.
 We consider the following elementary transport equation    \begin{equation}
 \label{e:transport}
 \left\{
 \begin{array}{ll}
 \d_t u - V u = 0, & (t, x) \in \R \times \M , \\
 u|_{t=0} = \varphi& x \in \M , 
 \end{array}
 \right.
 \end{equation}
 where $\M$ is a compact manifold, endowed with a density $|dx|$, and $V$ is a smooth vector field on $\M$.
A general question arising from physics concerns the large time behavior of the associated solution $u(t)=e^{tV}\varphi$ (where, with a slight abuse of notation, $(e^{tV})_{t\in \R}$ denotes the evolution group generated by $V$, say on $L^p(\M)$, $p<\infty$).
Assuming that $\div_{|dx|}(V)=0$, then one notices that~\eqref{e:transport} is conservative, e.g. in the sense that $\| u(t) \|_{L^p(\M)}=  \| \varphi \|_{L^p(\M)}$ for any $p \in [1,\infty]$. 
 Hence there is no hope to see convergence or dispersion of solutions to~\eqref{e:transport} in strong topologies.
One may rather consider these phenomena in a weak sense, and investigate, for instance for $\varphi \in C^\infty(\M)$ convergence in the sense of distributions: 
$$
 u(t) \rightharpoonup P_0\varphi \quad \text{ in } \mathcal{D}'(\M) , 
$$
where the equilibrium state $P_0\varphi$ is to be determined (if it exists).
This question is of particular interest in dynamical systems in which case $\M = SM$ is the unit sphere bundle over a compact Riemannian manifold $M$, and $V$ is the geodesic vector field acting on $\M=SM$.

This  last question has received attention in the case where $M$ has negative sectional curvature. In this case, the flow is ergodic with respect to the Liouville measure $L_g$, the only invariant functions are constant on $\M$, and one has $$P_0\varphi= \frac{1}{L_g(\M)} \int_\M \varphi dL_g .$$
It is by now well-understood (Ratner~\cite{Ratner87}, Dolgopyat~\cite{Dolgopyat98}, Liverani~\cite{Liverani04}, Tsujii~\cite{Tsujii10}, Nonnenmacher-Zworski~\cite{NonnenmacherZworski15}...) that this dynamical system is exponentially mixing, that is to say, that 
$$ \left| \left< u(t) -P_0\varphi , \psi \right>\right| \leq Ce^{-\eps t} \| \varphi\|_{\mathcal{H}}\|\psi\|_{\mathcal{H}'},$$
for test functions $\varphi,\psi$ taken in appropriate spaces $\mathcal{H}, \mathcal{H}'$ of distributions on $\M$.

Our last question is the counterpart to the above question on the torus, that is to say if $M = \T^d$, $\M=S\T^d=\T^d\times \IS^{d-1} \ni (x,\theta)$, and the geodesic vector field is simply given by $V= \theta\cdot \d_x$.

\begin{ques} 
\label{Q3}
On the torus, what is the asymptotics of the correlation function: for $\varphi, \psi \in C^\infty(S\T^d)$: does 
$$
\Cor_{\varphi,\psi}(t) : = \left<e^{tV}\varphi,\psi \right> 
$$
converge as $t \to +\infty$? What is its limit? 
What is $P_0\varphi$? What is the convergence rate? Does one have an asymptotic expansion? What are the (possibly anisotropic) spaces $\mathcal{H}$ adapted to the convergence?  
  \end{ques}
 Note that when the phase space is $T\T^d = \IT^d\times\IR^d$ instead of $S\T^d=\T^d\times \mathbb{S}^{d-1}$, this is a classical question studied in kinetic theory, see e.g.~\cite{MV11} and the references therein.
In that setting, convergence to zero holds at rate $t^{-\infty}$ (and even $e^{-\gamma t}$ if $\varphi, \psi$ are taken in appropriate spaces of analytic functions). Here, we shall see that the restriction to the energy shell $\IS^{d-1}$ (instead of the whole $\R^d$) is responsible for a different asymptotic behaviour.

As we shall see in the sequel, Question~\ref{Q3} is actually much simpler than the above two geometric questions. However, a fine answer to Question~\ref{Q3} is the cornerstone to takle Questions~\ref{Q1} and~\ref{Q2}.

 %%%%%%%%%%%%%%%%%%%%%%%%%%%%%
\section{Three results}
In this section, we briefly present some of the results we obtain concerning the above Questions~\ref{Q1},~\ref{Q2} and~\ref{Q3}.
\subsection{Geometric zeta function associated to orthogeodesics on the torus}
\label{s:zeta}
As a first result, we obtain the following proposition.
\begin{proposition}
\label{l:starting-point}
There is $T_0>0$ large enough, such for any $T>T_0$ the subset 
$$
\left\{\gamma\in\mathcal{P}_{K_1,K_2}: T_0<\ell(\gamma) \leq T \right\}
$$ 
of $\mathcal{P}_{K_1,K_2}$ is finite and we have as $T \to + \infty$
\begin{equation}\label{e:asymptotic-counting}
 \sharp \left\{\gamma\in\mathcal{P}_{K_1,K_2}: T_0<\ell(\gamma)\leq T\right\}=\frac{\pi^{\frac{d}{2}}T^d}{(2\pi)^{d}\Gamma\left(\frac{d}{2}+1\right)}+\mathcal{O}(T^{d-1}).
\end{equation}
\end{proposition}
In the case where $K_1=K_2=\{0\}$, this exactly amounts to count the number of lattice points in $2\pi\IZ^d$ of length less than $T$ and understanding the optimal size of the remainder in~\eqref{e:asymptotic-counting} is a deep problem in number theory. As a particular consequence of (a weak form of)~\eqref{e:asymptotic-counting}, the series $s \mapsto \zeta_{\Conv}(K_1,K_2,s)$ is absolutely convergent and holomorphic in $\{\Re(s)>d\}$.
The next step is to describe its meromorphic continuation to $\mathbb{C}$:
\begin{theorem}\label{t:maintheo-Mellin} Let $K_1$ and $K_2$ be two strictly convex and compact subsets of $\IR^d$ ($d\geq 2$) with smooth boundary, 
 then
 $$s\in \{\operatorname{Re}(s)>d\}\mapsto \zeta_{\Conv}(K_1,K_2,s)$$
 extends meromorphically to $\mathbb{C}$, its poles are located at $s=1,\ldots,d$ and are simple.
\end{theorem}

As for classical zeta functions in number theory, it is natural to compute the explicit values of the residues and, due to the geometric nature of the problem, one would like to express them in terms of natural geometric quantities attached to the convex sets $K_1$ and $K_2$. As a preliminary, we recall Steiner's formula~\cite[\S4]{Schneider14}. Denote by $B_d$ the unit ball in $\R^d$. For a compact and convex subset $K$ of $\IR^d$, $t \mapsto \text{Vol}_{\IR^d}\left(K+tB_d\right)$ is a polynomial of degree $d$ with nonnegative coefficients. The formula 
\begin{equation}\label{e:steiner}\text{ for all } t>0,\quad \text{Vol}_{\IR^d}\left(K+tB_d\right) =\sum_{\ell=0}^dV_{d-\ell}\left(K\right)\frac{\pi^{\frac{\ell}{2}}}{\Gamma\left(\frac{\ell}{2}+1\right)} t^{\ell},\end{equation}
thus defines the family of coefficients $V_{\ell}\left(K\right)\geq0$, called the $\ell$-\emph{intrinsic volume} of the convex set $K$. 
Note that $V_0(K)=1$, $V_d(K)=\text{Vol}_{\IR^d}(K)$. Moreover, if $\partial K$ has smooth boundary, one finds by the Minkowski-Steiner formula~\cite[\S4.2]{Schneider14}~\cite[p.~86]{teissier2016bonnesen}:
$$V_{d-1}(K)=\frac{1}{2}\text{Vol}(\partial K),$$
where $\text{Vol}$ is the $(d-1)$-volume measure on $\partial K$ induced by the Euclidean structure on $\IR^d$. Observe that $V_{d-\ell}\left(K\right)=0$ for any $0\leq \ell\leq d-1$ when $K$ is reduced to a point. Other properties of these intrinsic volumes are their invariance under Euclidean isometries (i.e. any composition of a rotation with a translation), their continuity with respect to the Hausdorff metric and their additivity\footnote{A functional satisfying such an additive property is referred as a \emph{valuation}~\cite[\S6.1]{Schneider14}.} on convex subsets of $\IR^d$, i.e.
$$\forall\ 0\leq \ell \leq d,\quad V_{\ell}\left(K\right)+V_{\ell}\left(L\right)=V_{\ell}\left(K\cup L\right)+V_{\ell}\left(K\cap L\right),$$
whenever $K$, $L$, $K\cup L$, $K\cap L$ are convex subsets of $\IR^d$. In fact, a classical Theorem of Hadwiger states that any functional on the convex subsets of $\IR^d$ enjoying these three properties is a linear combination of these intrinsic volumes~\cite[Th.~6.4.14]{Schneider14}.

Our second main theorem expresses the residues of $\zeta_{\Conv}(K_1,K_2,s)$ in terms of these intrinsic volumes:
\begin{theorem}\label{t:maintheo-residues}
 Let $K_1$ and $K_2$ be two strictly convex and compact subsets of $\IR^d$ ($d\geq 2$) with smooth boundary. Suppose in addition that $\Sigma_1=\mathfrak{p}(\partial K_1)$ (resp. $\Sigma_2=\mathfrak{p}(\d K_2)$) is oriented by the outgoing (resp. ingoing) normal vector to $K_1$ (resp. $K_2$).
 Then, the function
 $$s\mapsto \zeta_{\Conv}(K_1,K_2,s)-\frac{1}{(2\pi)^d}\sum_{\ell=1}^d\frac{\pi^{\frac{\ell}{2}}\ell}{\Gamma\left(\frac{\ell}{2}+1\right)}\frac{V_{d-\ell}\left(K_1-K_2\right)}{s-\ell}$$
 extends holomorphically from $\{\operatorname{Re}(s)>d\}$ to $\mathbb{C}.$
\end{theorem}
Note that $-K_2$ is a convex set and thus so is $K_1-K_2$.
In the case where both $K_1$ and $K_2$ are reduced to points, this theorem recovers a classical property of Epstein zeta functions~\cite{Epstein03} as all mixed volumes vanish except for $V_0$. 

%%%%%%%%%%%%%%%%%%%%%%%%%%%%%%%%%%%%%%%
\subsection{Geometric Poisson formul{\ae} associated to orthogeodesics on the torus}
\label{s:poisson-ortho}

As far as the ortholength distribution  $\mathsf{T}_{K_1,K_2}$ is concerned, our first result describes the distributional singularities of $\widehat{\mathsf{T}}_{K_1,K_2}$.
\begin{theorem}[Poisson type formula]\label{e:maintheo-realaxis} Let $K_1$ and $K_2$ be two strictly convex and compact subsets of $\IR^d$ ($d\geq 2$) with smooth boundary. 
Then, with $\mathsf{T}_{K_1,K_2}$ defined in~\eqref{d:TKK}, we have 
$$
\singsupp \widehat{\mathsf{T}}_{K_1,K_2} \subset \Sp(-\sqrt{-\Delta} ) \cup \Sp(\sqrt{-\Delta} ),
$$
where $\Delta = \sum_{j=1}^d \d_{x_j}^2$ is the flat Laplace operator on $\T^d$.
\end{theorem}
 Here, we recall that the singular support of a distribution $\mathsf{T}$ is the complementary of the open set where the distribution is $C^\infty$. In particular, the singular support of the geometric distribution $\widehat{\mathsf{T}}_{\beta,K_1,K_2}$ is given by the eigenvalues of the Laplacian and it does not depend on the convex sets $K_1,K_2$. 
  A more precise form of the theorem is actually proven in~\cite{DangLeautaudRiviere21}, in which the form of the singularities is described (depending on the dimension $d$ only), and the leading coefficients explicited (which depend on the convex sets $K_1,K_2$). In particular, the mixed volumes involved in Theorem~\ref{t:maintheo-residues} appear here in the singularity at zero.

 Theorem~\ref{e:maintheo-realaxis} has a similar flavour as the wave trace formula except that the correspondence is somehow in the other sense and that it involves orthogeodesics of two given convex sets. More precisely, we start from the length orthospectrum between two convex sets, we then form the series $\widehat{\mathsf{T}}_{K_1,K_2}(t)=\sum_\gamma  e^{-it\ell(\gamma)}$, and its singular support coincides with the quantum spectrum $\Sp(\pm \sqrt{-\Delta})$ where $\Delta$ is the Laplacian. Another notable difference is that the singularities are more complicated in the sense that they involve distributions of the form $(t\pm \lambda-i 0)^{-k}$ with $k\geq 1$ that may not even be an integer if $d$ is even. We emphasize that, as in the Chazarain--Duistermaat--Guillemin formula, the singularities of $\widehat{\mathsf{T}}_{K_1,K_2}$ are {\em not} purely Dirac type distributions (and their derivatives). This is due to the fact that the counting measure $T_{K_1,K_2}$ is supported on the half--line, hence its Fourier transform $\widehat{\mathsf{T}}_{K_1,K_2}$ must have its ($C^\infty$ and analytic) wave front set contained in the half cotangent cone $\{(t;\tau); \tau<0\}\subset T^*\mathbb{R}$. This prevents the presence of purely $\delta^{(k)}(t)$--like singularities whose contribution to the wave front set would contain both positive and negative frequencies $\tau$. 

For all of these reasons, as far as Question~\ref{Q2} is concerned, the distribution $\widehat{\mathsf{T}}_{K_1,K_2}$ is not a Dirac comb (at all). We can partially remedy this issue as follows.
In fact, in view of having simple singularities and motivated by the recent developments on crystalline measures~\cite{Meyer2022}, one can
\begin{itemize}
\item twist the definition of $\mathsf{T}_{K_1,K_2}(t)$ with a weight $e^{i\int_\gamma\beta}$, where 
\begin{align}
\label{e:def-beta}
\beta = \sum_{j=1}^d \beta_j dx_j + d f , \quad \beta_j \in \R, \quad f \in C^\infty(\T^d;\R)
\end{align} is a (particular) closed one-form, in order to erase the singularity at zero; 
\item  symmetrize and renormalize the distribution $\mathsf{T}_{K_1,K_2}(t)$ in order to recover Dirac-type singularities.
\end{itemize}
This is the content of our last main result which extends in our geometric setup the Guinand--Meyer summation formula~\cite[Th.~5]{Meyer2016}.

 \begin{theorem}[Guinand--Meyer type formula]\label{e:GuinandWeil} Let $K_1$ and $K_2$ be two strictly convex and compact subsets of $\IR^d$ ($d\geq 2$) with smooth boundary and let $\beta$ as in~\eqref{e:def-beta} with $(\beta_1,\dots, \beta_d ) \notin \Z^d$. Let $\mu$ be the complex measure defined as 
\begin{eqnarray*}
\mu(t)= \sum_{\gamma\in\mathcal{P}_{K_1,K_2}:\ell(\gamma)>T_0}\frac{e^{i\int_\gamma\beta}}{\ell(\gamma)^{\frac{d-1}{2}}} \delta_{\ell(\gamma)}+(-i)^{d-1}\sum_{\gamma\in\mathcal{P}_{K_2,K_1}:\ell(\gamma)>T_0}\frac{e^{-i\int_\gamma\beta}}{\ell(\gamma)^{\frac{d-1}{2}}} \delta_{-\ell(\gamma)},
\end{eqnarray*} 
where we take the same orientation conventions for\footnote{In particular, both sets are a priori distinct.} $\mathcal{P}_{K_2,K_1}$ and $\mathcal{P}_{K_1,K_2}$. 

Then, there exist complex numbers $(c_\lambda)_{\lambda\in \operatorname{Sp}(\pm \sqrt{-\Delta_{[\beta]}})}$ and $r$ belonging to $L^p_{\operatorname{loc}}(\IR)$ for every $1\leq p<\infty$ such that
\begin{eqnarray*}
 \widehat{\mu}(\tau)=\sum_{\lambda\in \operatorname{Sp}(\pm \sqrt{-\Delta_{[\beta]}})} c_\lambda \delta_\lambda+r.
\end{eqnarray*} 
\end{theorem}
In the case $(\beta_1,\dots, \beta_d ) \in \Z^d$, the result would be similar except for an extra singularity at $\tau=0$ that may be more singular than the Dirac distribution. 
In the case where $K_1$ and $K_2$ are distinct points and where $d=3$, it was in fact proved that $r\equiv 0$ in~\cite[Th.~5]{Meyer2016}. We also recover through our geometric setup the result of~\cite[\S2]{LevReti2021} concerning the construction of  \emph{crystalline distributions} in higher dimensions (when $d$ is odd). We finally prove that $r$ is not identically $0$ as soon as $d\geq 5$, even in the case where $K_1,K_2$ are points. See also~\cite{Guinand59,LevOlevskii2016,RadchenkoViazovska2019} for earlier related results and~\cite{Meyer2022} for a review on recent developments in that direction.

%%%%%%%%%%%%%%%%%%%%%%%%%%%%%%%%%%%%%%%
\subsection{Decay of the correlation function and emergence of quantum dynamics}
\label{s:emergence-quantum}

Theorems~\ref{t:maintheo-mellin-function} and~\ref{t:maintheo-resolvent-function} (as well as their analogues in the case of differential forms) are consequences of the fact that we can give a full expansion of the Schwartz kernel of the geodesic flow. 
For simplicity, we may interpret the correlation function $\Cor_{\varphi,\psi}(t)$ in the sense of the space $L^2(S\T^d)$ for the Liouville measure, which is simply $dL_g = dx \otimes d\Vol(\theta)$ on $S\T^d=\T^d\times \IS^{d-1}$.
For instance, the first term in the asymptotic expansion reads as follows.  
  
\begin{theorem}[Time asymptotics of the geodesic flow, function case, leading term]
\label{t:maintheo-correlations} 
For every smooth functions $\varphi, \psi\in\ml{C}^{\infty}(S\IT^d)$, we have 
 \begin{align*}
  \Cor_{\varphi,\psi}(t) &= \left(P_0 \varphi, P_0 \psi \right)_{L^2(\IS^{d-1})} + \left( \frac{2\pi}{t}\right)^{\frac{d-1}{2}}  \sum_{\pm}  e^{\mp i\frac{\pi}{4}(d-1)} \left( \frac{e^{\pm i t\sqrt{-\Delta}}}{\sqrt{-\Delta}^{\frac{d-1}{2}}} \Pi_0^\pm \varphi, \Pi_0^\pm \psi \right)_{L^2(\T^d)} \\ 
  &  \quad+ \mathcal{O}_{\varphi, \psi} \left(  \frac{1}{t^{\frac{d-1}{2}+1}} \right)
\end{align*}
  where  $(P_0 f)(\theta) = \frac{1}{(2\pi)^d}\int_{\T^d}f(x,\theta) dx$ is the projection onto the vector space of $x$-invariant functions and
   $$ 
   (\Pi_0^\pm f)(x) =  \sum_{\xi\neq0} \hat{f}_\xi \left( \pm \frac{\xi}{|\xi|} \right) e^{i \xi \cdot x},\quad\text{for}\quad f(x,\theta) =  \sum_{\xi\in \Z^d} \hat{f}_\xi \left(\theta \right) e^{i \xi \cdot x}. 
    $$
    
Equivalently, for every smooth function $\psi\in\ml{C}^{\infty}(S\IT^d)$, one has
\begin{align*}
t^{\frac{d-1}{2}}\left(\psi\circ e^{-tV}(x,\theta)-\frac{1}{(2\pi)^d}\int_{\IT^d}\psi(y,\theta)dy\right)&= (2\pi)^{\frac{d-1}{2}}\sum_{\pm}\mathbf{P}_{\pm}^\dagger\frac{e^{\pm i\left(t\sqrt{-\Delta}-\frac{\pi}{4}(d-1)\right)}}{(-\Delta)^{\frac{d-1}{4}}}\mathbf{P}_{\pm}\\
&\quad +\ml{O}_{\ml{D}^{\prime}(S\IT^d)}(t^{-1})
\end{align*}
where $\Delta=\sum_{j=1}^d\partial_{x_j}^2$ is the Euclidean Laplacian on $\IT^d$, 
$$\mathbf{P}_{\pm}:\psi\in\ml{C}^{\infty}(S\IT^d)\mapsto\sum_{\xi\neq 0}\frac{1}{(2\pi)^d}\int_{\IT^d}\psi\left(y,\pm\frac{\xi}{|\xi|}\right)e^{i(y-x)\cdot\xi}dy\in\ml{C}^{\infty}(\IT^d)$$
and
 $$\mathbf{P}_{\pm}^\dagger:f\in\ml{C}^{\infty}(\IT^d)\mapsto \sum_{\xi\neq 0}\frac{1}{(2\pi)^d}\left(\int_{\IT^d}f\left(y\right)e^{i(y-x)\cdot\xi}dy\right)\delta_0\left(\theta\mp\frac{\xi}{|\xi|}\right)\in\ml{D}^{\prime}(S\IT^d)$$
\end{theorem}
We actually provide with a full asymptotic expansion and a precise description of the remainder terms at each step. This result could (and will) in fact be expressed in terms of anisotropic Sobolev norms. 

In order to keep track of the comparison with negatively curved manifolds, such a result can be viewed as a simple occurence of the emergence of quantum dynamics (through the half-wave group $(e^{\pm it \sqrt{-\Delta}})_{t \in \R}$ on the torus) in the long time dynamics of geodesic flows (i.e. $(e^{tV})_{t\in\R}$ on $S\T^d$).
This might be interpreted as the fact that the classical dynamics, corrections to the leading order in the convergence towards the equilibrium are governed by quantum fluctuations (with amplitude decaying in time like $t^{-\frac{d-1}{2}}$).

 This phenomenon was recently exhibited by Faure and Tsujii in the general context of contact Anosov flows~\cite{FaureTsujii15, FaureTsujii17, FaureTsujii17b, FaureTsujii21}. See also~\cite{DyatlovFaureGuillarmou2015} for related results of Dyatlov, Faure and Guillarmou in the particular case of geodesic flows on \emph{hyperbolic} manifolds. Compared with the results of Faure and Tsujii, we emphasize that our analysis heavily relies on the algebraic structure of our flows as in the hyperbolic settings treated in~\cite{Ratner87, DyatlovFaureGuillarmou2015}. Moreover, we are dealing with completely integrable systems which have in some sense opposite behaviours compared with the dynamical situations considered in all these references. In particular, due to the integrable nature of our system, the asymptotic expansion in terms of the quantum propagator is polynomial rather than exponential as in~\cite[Th.~1.2]{FaureTsujii21}. This is reminiscent of the much weaker mixing properties of the geodesic flow in this situation.

%%%%%%%%%%%%%%%%%%%%%%%%%%%%%%%
\section{Three remarks on the proofs}
In this section, we briefly discuss some points of the proofs, with an  emphasis on their analytic side (and completely omitting their geometric side, for which we refer to~\cite{DangLeautaudRiviere21}). 

\subsection{Decay of the correlation function}
In this section, we briefly sketch the proof of Theorem~\ref{t:maintheo-correlations}, which only involves classical stationary/non-stationary phase estimates, and, at the same time, is at the root, the heart  and the cornerstone of the proofs of Theorems~\ref{t:maintheo-Mellin}, \ref{t:maintheo-residues} and \ref{e:GuinandWeil}.

To give an asymptotic expansion of the correlation function $\Cor_{\varphi,\psi}(t)$, we expand in partial Fourier series in $x \in \T^d$ the functions 
 $$\psi(x,\theta)=\sum_{\xi\in\IZ^d}\widehat{\psi}_\xi(\theta)e^{i\xi \cdot x},  \quad \varphi(x,\theta)=\sum_{\xi\in\IZ^d}\widehat{\varphi}_\xi(\theta)e^{i\xi \cdot x} .$$
 Recalling that the geodesic flow is simply the shift
 $$
e^{tV} :  (x,\theta) \mapsto (x+t\theta, \theta) ,
 $$
 and using the Plancherel identity, we thus have 
\begin{align*}
\Cor_{\varphi,\psi}(t)& = \left(e^{tV}\varphi,\psi \right)_{L^2(S\T^d)}= \int_{\T^d \times \IS^{d-1}}\varphi (x +t\theta , \theta)\overline{\psi}(x,\theta)  dx \, d\Vol(\theta) \\
 & =\sum_{\xi\in\IZ^d}\int_{\IS^{d-1}} e^{it \xi\cdot\theta}
\widehat{\varphi}_{\xi}(\theta)\overline{\widehat{\psi}_{\xi}}(\theta) d\Vol(\theta) \\
 & =\underbrace{\int_{\IS^{d-1}}  
\widehat{\varphi}_{0}(\theta)\overline{\widehat{\psi}_{0}}(\theta)  d\Vol(\theta)}_{\quad =\left(P_0 \varphi, P_0 \psi \right)_{L^2(\IS^{d-1})}} +\sum_{\xi\in\IZ^d, \xi \neq 0 }\underbrace{\int_{\IS^{d-1}}e^{it \xi\cdot\theta} 
\widehat{\varphi}_{\xi}(\theta)\overline{\widehat{\psi}_{\xi}}(\theta)  d\Vol(\theta)}_{\text{oscillatory integral}}
\end{align*}
The time-invariant term is already identified and it remains to give an asymptotic expansion of the oscillatory integrals, which rewrite under the form
$$I_F(\xi,t):=\int_{\IS^{d-1}}e^{it\xi\cdot\theta}F(\theta) d\Vol(\theta),\quad \text{as }t\to+\infty . $$
Up to a rotation of the vector $\xi$, we may assume that $\xi=|\xi|e_d$ where $e_d$ is the last vector of the canonical basis of $\R^d$ and  $\lambda =t |\xi|$ is large (recall that $|\xi|\geq 1$ in the sum), so that we are left to study 
 $$I_F(\xi,t):=\int_{\IS^{d-1}}e^{i\lambda e_d \cdot\theta}F(\theta)d\Vol(\theta),\quad\text{as }\lambda\to+\infty. $$
 Estimating these kind of integrals is a classical topic in harmonic analysis, see e.g.~\cite{Herz62, Littman63} for a rough estimate, and~\cite[Th.~7.7.14]{Hormander90}, \cite[Section~VIII-3, p347]{Steinbook} and~\cite[Th.~3.38, p140]{DyatlovZworski19} for fine asymptotic expansions. 
In this oscillatory integral, the phase is the height function $\IS^{d-1} \ni \theta \mapsto e_d \cdot \theta$. The latter has only two (nondegenerate) critical points  given by $\theta=\pm e_d$.

Nonstationary  phase expansion yields (slightly informally) ``away from the poles $\pm e_d$'':
 \begin{align}
 \label{e:away-poles}
 |I_F(\xi,t)| \leq \frac{C_N}{\lambda^N} \|F\|_{W^{N,1}} =  \frac{C_N}{t^N}  \frac{1}{|\xi|^N} \|\widehat{\varphi}_{\xi} \overline{\widehat{\psi}_{\xi}}\|_{W^{N,1}(\IS^{d-1})} \leq  \frac{C_N}{t^N}  \frac{1}{|\xi|^N} \|\widehat{\varphi}_{\xi}\|_{H^N(\IS^{d-1})} \| \widehat{\psi}_{\xi} \|_{H^N(\IS^{d-1})} .
 \end{align}
Near the poles $\pm e_d$, a stationary  phase expansion (we omit the details here, see~\cite{DangLeautaudRiviere21}) yields, for differential operators $L_{j}^\pm$ of order $2j$ on $\IS^{d-1}$ (with $L_{0}^\pm=$ identity):
\begin{align*}
I_F(\xi,t) & = 
  e^{\pm i \lambda} e^{\mp i\frac{\pi}{4}(d-1)}\left(\frac{2\pi}{\lambda}\right)^{\frac{d-1}{2}}
  \sum_{j=0}^{N-1} \frac{1}{\lambda^j} L_{j}^\pm \left(F\right) \left(\pm e_d \right) + \ml{O}_N\left( \frac{1}{\lambda^{N+\frac{d-1}{2}}} \right) \left\| F\right\|_{W^{2N+ d,1}} . 
\end{align*}
  Back to the correlation function, and recalling that we made a rotation so that $e_d = \frac{\xi}{|\xi|}$, we have thus obtained 
 \begin{align}
 \label{e:poles}
\Cor_{\varphi,\psi}(t)  
 & =  \left(P_0\varphi,P_0\psi \right)_{L^2} + \sum_{\xi\neq 0}
e^{\pm it|\xi|} e^{\mp i\frac{\pi}{4}(d-1)}\left(\frac{2\pi}{t|\xi|}\right)^{\frac{d-1}{2}}
  \sum_{j=0}^{N-1} \frac{1}{(t|\xi|)^j} L_{j, \frac{\xi}{|\xi|}}^\pm \left(\widehat{\varphi}_{\xi} \overline{\widehat{\psi}_{\xi}}\right) \left(\pm \frac{\xi}{|\xi|}\right)  \nonumber \\
& \quad  + \ml{O}_N\left(t^{-N-\frac{d-1}{2}}\right)  \sum_{\xi\neq 0}\frac{\left\|\widehat{\varphi}_{\xi}\right\|_{H^{2N+ d}(\IS^{d-1})}\left\| \widehat{\psi}_{\xi}\right\|_{H^{2N+ d}(\IS^{d-1})}}{|\xi|^{N+\frac{d-1}{2}}} .
\end{align}
This concludes the sketch of the proof of a higher order version of the results of Theorem~\ref{t:maintheo-correlations}.

%%%%%%%%%%%%%%%%%%%%%%%%%%%%%%%%%%%%%%%%
\subsection{Anisotropic Sobolev spaces, Mellin, Laplace, and Fourier transforms}
\label{s:sobobo}
When observing the main terms and the remainder in the above asymptotic expansions~\eqref{e:away-poles} and~\eqref{e:poles}, and recalling that $\xi \in \Z^d$ is the Fourier variable with respect to $x \in \T^d$, we see that different regularity/singularity is required/allowed in the directions of  $\theta \in \IS^{d-1}$ and $x \in \T^d$.
Also, different regularity is required for $\theta$ near and away from the poles $\pm\frac{\xi}{|\xi|}$, but we shall not discuss this issue here.
This leads us to define anisotropic Sobolev spaces of distributions on $S\IT^d$ as follows
$$\mathcal{H}^{M,N}(S\IT^d):=\left\{u\in\ml{D}^\prime(S\IT^d):\ \|u\|_{\mathcal{H}^{M,N}(S\IT^d)}^2 := \sum_{\xi\in\IZ^d} \langle \xi \rangle^{2N}\|\widehat{u}_\xi\|_{H^M(\IS^{d-1})}^2<+\infty\right\}, $$
where $(M,N)\in\IR^2$, $\langle \xi \rangle = (1+|\xi|^2)^{\frac{1}{2}}$, and 
$$u(x,\theta)=\sum_{\xi\in\IZ^d}\widehat{u}_\xi(\theta)\frac{e^{i\xi\cdot x}}{(2\pi)^{\frac{d}{2}}},$$
with $\widehat{u}_\xi\in\ml{D}^{\prime}(\IS^{d-1})$, and where $\|.\|_{H^{M}}$ denotes the standard Sobolev norm on $\IS^{d-1}$. Roughly speaking, $u=u (x,\theta)\in \mathcal{H}^{M,N}(S\IT^d)$ if $u$ has $H^N$ regularity in the variable $x \in \T^d$ (which has a geometric interpretation as the so-called horizontal direction) and $H^{M}$ regularity in the variable $\theta \in \IS^{d-1}$  (which has a geometric interpretation as the so-called vertical direction).
With this notation at hand, the remainder in~\eqref{e:poles} may for instance be rewritten as 
$$
 \sum_{\xi\neq 0}\frac{\left\|\widehat{\varphi}_{\xi}\right\|_{H^{2N+ d}}\left\| \overline{\widehat{\psi}_{\xi}}\right\|_{H^{2N+ d}}}{|\xi|^{N+\frac{d-1}{2}}}  \leq  \left\| \varphi \right\|_{\H^{2N+ d,-(N+\frac{d-1}{2})}}
 \left\| \psi \right\|_{\H^{2N+ d,-(N+\frac{d-1}{2})}} , 
$$
and we notice that the expansion~\eqref{e:poles}, up to order $N$ holds for distributions $\varphi,\psi$ being only $H^{-\frac{d-1}{2}}$ in $x \in \T^d$ but having regularity $H^{2N+d}$ in $\theta \in \IS^{d-1}$. Note that the restriction to $H^{-\frac{d-1}{2}}$ along the $x$-variable comes from the leading term in the asymptotic expansion.

\bigskip
Motivated by and in analogy with Theorems~\ref{t:maintheo-Mellin} and~\ref{e:maintheo-realaxis}, given two functions $\varphi,\psi \in C^\infty(\T^d\times \IS^{d-1})$, we now discuss properties of the Mellin, Laplace and Fourier transforms of the associated correlation function $\Cor_{\varphi,\psi}$, namely
\begin{align*}\M\big(\Cor_{\varphi,\psi}\big) (s) &:= \int_1^\infty t^{-s} \, \Cor_{\varphi,\psi}(t) dt \quad \text{(Mellin transform)}, \\
\ml{L}\big( \Cor_{\varphi,\psi}\big)(s) &:= \int_0^\infty e^{-st} \, \Cor_{\varphi,\psi}(t) dt \quad \text{(Laplace transform)} ,\\
 \F\big(\Cor_{\varphi,\psi}\big) (\tau) &:= \int_0^\infty e^{-i\tau t} \, \Cor_{\varphi,\psi}(t) dt \quad \text{(Fourier transform)} .
 \end{align*}
More precisely, we may consider them as maps $\mathcal{C}^{\infty}(S\IT^d)\rightarrow \mathcal{D}^{\prime}(S\IT^d)$, applied to $\varphi$, and tested against $\psi$, namely e.g. $\M\big(\Cor_{\varphi,\psi}\big) (s) =\left< \M(s)\varphi , \psi \right>$.
Concerning the Mellin transform, our first result reads as follows.

\begin{theorem}[Mellin transform, function case]\label{t:maintheo-mellin-function} Let $\chi\in\ml{C}^\infty_c([1,+\infty))$ such that $\chi=1$ in a neighborhood of $1$ and let $N\in\IZ_+$. Then, the operator
$$\mathcal{M}(s):=\int_1^{\infty}t^{-s}e^{-tV*}|dt|:\mathcal{C}^{\infty}(S\IT^d)\rightarrow \mathcal{D}^{\prime}(S\IT^d)$$
splits as
$$\mathcal{M}(s)=\mathcal{M}_0(s)+\mathcal{M}_\infty(s),$$
where
$$\mathcal{M}_0(s):=\int_1^{\infty}\chi(t)t^{-s}e^{-tV*}|dt|:\ml{H}^{N,-N/2}(S\IT^d)\rightarrow\ml{H}^{N,-N/2}(S\IT^d)$$
is a holomorphic family of bounded operators on $\mathbb{C}$ and where
$$\mathcal{M}_\infty(s):=\int_1^{\infty}(1-\chi(t))t^{-s}e^{-tV*}|dt|:\ml{H}^{N,-N/2}(S\IT^d)\rightarrow\ml{H}^{-N,N/2}(S\IT^d)$$
extends as a meromorphic family of bounded operators from $\{\operatorname{Re}(s)>1\}$ to $\{\operatorname{Re}(s)>1-N\}$
 with only a simple pole at $s=1$ whose residue is given by
 $$\forall\psi\in\mathcal{C}^{\infty}(S\IT^d),\quad\Res_{s=1}\left(\mathcal{M}_\infty(s)\right)(\psi)(x,\theta)=\frac{1}{(2\pi)^d}\int_{\IT^d}\psi(y,\theta)dy.$$
\end{theorem}
In particular, this result tells us that the operator $\mathcal{M}(s):\mathcal{C}^{\infty}(S\IT^d)\rightarrow \mathcal{D}^{\prime}(S\IT^d)$
extends meromorphically from $\{\operatorname{Re}(s)>1\}$ to the whole complex plane with only a simple pole at $s=1$. Yet, the statement is more precise as it allows us to describe the allowed regularity for this meromorphic continuation. We emphasize that the mapping properties of $\mathcal{M}_0(s)$ are rather immediate from the definition of our anisotropic norms and the main difficulty in this statement is about the ``regularizing'' properties of $\mathcal{M}_\infty(s)$ along the $x$-variable that will be instrumental in our applications to convex geometry. The proof of this result essentially relies on the asymptotic expansion~\eqref{e:away-poles}, together with the microlocal smoothing effect in the direction of the flow coming from the integration over time in $\M(s)$.
Remarkably enough, the meromorphic continuations are valid on spaces of distributions that are regular along the vertical bundle to $S\IT^d$ (i.e. the tangent space to $\IS^{d-1}$) and that may have negative Sobolev regularity along the horizontal bundle (i.e. the tangent space to $\IT^{d}$). Moreover, the output of $\mathcal{M}_\infty(s)$ will have some regularity along the $x$-variable. In particular, the anisotropic Sobolev spaces $\mathcal{H}^{N,-N/2}$ contain the Dirac distribution $\delta_{[0]}(x)$ for $N>d$, and this is typically the kind of distributions that we will pick as test functions in order to derive our main applications in convex geometry (see~\eqref{e:current-zeta} below). In order to prove Theorems~\ref{t:maintheo-Mellin} and~\ref{t:maintheo-residues}, we will in fact need to prove more general statements for the action of $\mathcal{M}(s)$ on differential forms or more precisely on certain anisotropic Sobolev spaces of currents. Among other things, the action on differential forms will be responsible for the presence of the extra poles at $s=2,\ldots,d$ but this simplified statement already illustrates the kind of properties we are aiming at. 

The same spaces will also allow us to prove the following counterpart statements concerning the Laplace transform.
\begin{theorem}
\label{t:laplap}
The operator 
\begin{align}
\label{e:def-laplap}
\mathcal{L}(s):=(V+s)^{-1}=\int_0^{\infty}e^{-st}e^{-tV*}|dt|:\mathcal{C}^{\infty}(S\IT^d)\rightarrow \mathcal{D}^{\prime}(S\IT^d) ,
\end{align}
defined for $\Re(s)>0$, continues as a $\mathcal{C}^\infty$ function to
$$
\{ \Re(s)\geq 0 \}  \setminus \left( i\Sp(\sqrt{-\Delta}) \cup  i\Sp( - \sqrt{-\Delta})  \right) .
$$
\end{theorem}
 Taking the limit on the imaginary axis $s \to i \tau \in i\R$, $\Re(s)>0$, one can check that 
$$\ml{L}\big(\Cor_{\varphi,\psi}\big) (i \tau + \alpha) \to_{\alpha \to 0, \alpha>0} \F\big(\Cor_{\varphi,\psi}\big)(\tau)\  \text{in } \ml{S}'(\R),$$ 
and we directly obtain the following corollary.
\begin{corollary}
For all $\varphi,\psi \in C^\infty(S\T^d)$,  $\singsupp \F\big(\Cor_{\varphi,\psi} \big) \subset \Sp(\sqrt{-\Delta}) \cup \Sp(-\sqrt{-\Delta})$.
\end{corollary}
Note that when considering the case of trigonometric polynomials in the $x$-variable, one can verify that the singular support is 
included in a finite part of the Laplace spectrum corresponding to the Fourier modes under consideration. This corollary is of course (at least formally) related to Theorem~\ref{e:maintheo-realaxis}.
In the rough result of Theorem~\ref{t:laplap}, we did not keep track of the anisotropic Sobolev spaces, for the $\ml{C}^k$ regularity of the extension depends on the $N$ in $\ml{H}^{N,-N/2}(S\IT^d)$.  If we only consider the $\ml{C}^0$ regularity of the extension, an example of a more precise statement in which we keep track of the anisotropic Sobolev spaces and the description of the singularities is given by the following result.
\begin{theorem}[Laplace transform, function case, continuous continuation]
\label{t:maintheo-resolvent-function} Let $\chi\in\ml{C}^\infty_c([0,+\infty))$ such that $\chi=1$ in a neighborhood of $0$ and let $N\in 2\IZ_+^*+d$. Then, the operator $\mathcal{L}(s)$ in~\eqref{e:def-laplap}
splits as
$$\mathcal{L}(s)=\mathcal{L}_0(s)+\mathcal{L}_\infty(s),$$
where
$$\mathcal{L}_0(s):=\int_0^{\infty}\chi(t)e^{-st}e^{-tV*}|dt|:\ml{H}^{N,-N}(S\IT^d)\rightarrow\ml{H}^{N,-N}(S\IT^d)$$
is a holomorphic family of bounded operators on $\mathbb{C}$ and where
$$\mathcal{L}_\infty(s):=\int_0^{\infty}(1-\chi(t))e^{-st}e^{-tV*}|dt|:\ml{H}^{N,-N/2}(S\IT^d)\rightarrow\ml{H}^{-N,N/2}(S\IT^d)$$
extends continuously from $\{\operatorname{Re}(s)>0\}$ to 
 \begin{enumerate}
  \item $\{\operatorname{Re}(s)\geq 0\}\setminus\{0\}$ if $d\geq 4$,
  \item $\{\operatorname{Re}(s)\geq 0\}\setminus\{\pm i|\xi|:\xi\in\IZ^d\}$ if $d=2,3$.
 \end{enumerate}
Moreover, in any dimension, one has, as $s\rightarrow 0^+$,
 $$(V+s)^{-1}(\psi)(x,\theta)=\frac{1}{(2\pi)^ds}\int_{\IT^d}\psi(y,\theta)dy+\mathcal{O}_{\ml{D}^\prime}(1),$$
and, when $d=2,3$, one has, as $s\rightarrow \pm i|\xi_0|$ (with $|\xi_0|\neq 0$),
\begin{multline*}(V+s)^{-1}(\psi)(x,\theta)\\
=\frac{e^{\mp i\pi\frac{d-1}{4}}g_d(s\mp i|\xi_0|)}{(2\pi)^{\frac{d+1}{2}}|\xi_0|^{\frac{d-1}{2}}}\sum_{\xi:|\xi|=|\xi_0|}e^{i\xi\cdot x}\delta_0\left(\theta\mp \frac{\xi}{|\xi|}\right)\int_{\IT^d}\psi\left(y,\pm\frac{\xi}{|\xi|}\right)e^{-i\xi\cdot y}dy+\mathcal{O}_{\ml{D}^\prime}(1),
\end{multline*}
where
$$g_2(z):=\frac{\sqrt{2\pi}}{\sqrt{z}},\quad\text{and}\quad g_3(z):=-\ln(z).$$
\end{theorem}
Again,  Theorems~\ref{t:laplap} and~\ref{t:maintheo-resolvent-function} are consequences of the precise asymptotic expansion of the correlation function in~\eqref{e:poles} (see also Theorem~\ref{t:maintheo-correlations}), together with the trivial but crucial fact that, when formally taking the Laplace transform of~\eqref{e:poles}, we have that $\int_1^\infty \frac{e^{t(\pm i|\xi|-s)}}{t^{\frac{d-1}{2}}} dt$, defined originally on $\Re(s)>0$, extends smoothly to $\Re(s)\geq 0$ except for a singularity (depending on the dimension $d$) at $s=\pm i|\xi|$. Again we emphasize that this theorem provides us a smoothing effect in the $x$-variable that will be intrumental to derive our Poisson formulae in convex geometry. 
Note that, away from the singularities (and if we do not care about this smoothing property), rather than appealing to nonstationary phase arguments, we could as well have applied Mourre's commutator method~\cite{AmreinBoutetGeorgescu96} to the family (indexed by $\xi \in \Z^d$) of multiplication operators $u\in L^2(\IS^{d-1})\mapsto (\xi\cdot\theta) u(\theta)\in L^2(\IS^{d-1})$ (and then sum over $\xi \in \Z^d$).

 \subsection{Back to orthospectra: linking geometry and dynamics}
In this paragraph, we very roughly describe the strategy to prove Proposition~\ref{l:starting-point}, Theorems~\ref{t:maintheo-Mellin}, \ref{t:maintheo-residues} and \ref{e:GuinandWeil}, and in particular the links between these geometric questions and the dynamical issues of Section~\ref{s:emergence-quantum}.
Since the seminal work of Margulis~\cite{Margulis69, Margulis04}, it is well understood that on negatively curved manifolds, it is convenient to lift geometric problems on the manifold to its unit cotangent bundle. For instance, properties of Poincar\'e series are related to the asymptotic properties of the geodesic flow, and more specifically to its mixing properties. In a recent work~\cite{DangRiviere20d}, the first and last authors of the present note formulated this relation using the theory of De Rham currents and we follow this approach in the case of flat tori (although the curvature vanishes everywhere). Let us explain this connection without being very precise on the sense of the various integrals. We denote by $N_\pm(K_i)$ the outward/inward unit normal bundle to $K_i$ inside $S\IT^d$, namely 
$$N_\pm (K_i):=\left\{(\mathfrak{p}(x),d\mathfrak{p}(x)\theta),  x\in \partial K_i,  \pm \theta\text{ directly orthogonal to } \partial K_i \ \text{at}\ x\right\},$$
where $\mathfrak{p}$ is defined in~\eqref{can-proj}.
Then, given any nice enough function $\chi(t)$, say in $\mathcal{C}^{\infty}_c(\IR_+^*)$ (if we are interested in counting functions such as in Proposition~\ref{l:starting-point}), $t^{-s}$ (if we are interested in zeta function or in Mellin transforms) or $e^{-st}$ (if we are interested in resolvents or Laplace or Fourier transforms),  one has
\begin{equation*} \sum_{\gamma\in\mathcal{P}_{K_1,K_2}}\chi(\ell(\gamma))= \langle \mathsf{T}_{K_1,K_2} , \chi \rangle ,
 \end{equation*}
where $\mathsf{T}_{K_1,K_2}$ is defined in~\eqref{d:TKK}.
Adapting~\cite{DangRiviere20d}, we prove that 
\begin{equation}\label{e:current-zeta}\sum_{\gamma\in\mathcal{P}_{K_1,K_2}}\chi(\ell(\gamma))=(-1)^{d-1}\int_{S\IT^d}[N_{\sigma_1}(K_1)]\wedge \int_{\IR}\chi(t)e^{-t\mathbf{V}}\iota_V([N_{\sigma_2}(K_2)])dt ,
 \end{equation}
where $[N_{\sigma_i}(K_i)]$ is the De Rham current of integration on $N_{\sigma_i}(K_i)$, where $\sigma_i\in\{\pm\}$ depends on our orientation convention on each convex and where
$$e^{tV}:(x,\theta)\in S\IT^d \rightarrow (x+t\theta, \theta)\in S\IT^d$$
is the geodesic flow.  For simplicity, we drop the dependence in $\sigma_i$ for the end of this discussion.
Formula~\eqref{e:current-zeta} derives from the observation that elements in $\mathcal{P}_{K_1,K_2}$ are in one-to-one correspondance with the geodesic orbits in $S\IT^d$ joining the two Legendrian submanifolds $N(K_1)$ and $N(K_2)$.

This (rigorous) formula relates the counting of intersection points between two submanifolds with the geodesic flow acting on currents and can be formally (however slightly unrigorously) rewritten as 
$$
\langle \mathsf{T}_{K_1,K_2} , \chi \rangle  =  
\left\langle \Cor_{[N(K_1)],\iota_V[N(K_2)]}(t) , \chi(t) \right\rangle ,
$$
that is to say, in the sense of $\ml{S}'(\R)$
$$
\mathsf{T}_{K_1,K_2}(t)   =  \Cor_{[N(K_1)],\iota_V[N(K_2)]}(t) ,
$$
where $\Cor_{[N(K_1)],\iota_V[N(K_2)]}(t)$ denotes the correlation distribution (which is not a function of $t$ anymore) of the two De Rham integration currents $[N(K_1)]$ and $\iota_V[N(K_2)]$.
  Pairing the correlation distribution with $\chi(t)=e^{-st}$ (and then taking the limit $s=i\tau+\alpha$ with $\alpha\to 0^+$ to obtain the Fourier transform as the boundary value on $i\R$ of the Laplace transform as in~\cite[Th.3.1.11]{Hormander90}) and with $\chi(t)=t^{-s}$ relates in particular Questions~\ref{Q1},~\ref{Q2} and~\ref{Q3} by 
\begin{align*}
 \underbrace{\ml{L}\left(\mathsf{T}_{K_1,K_2}\right)(s)}_{\text{Laplace transform, related to Question~\ref{Q2}}}  & =\langle  \mathsf{T}_{K_1,K_2}, e^{-st} \rangle = \left\langle \underbrace{\Cor_{[N(K_1)],\iota_V[N(K_2)]}(t)}_{\text{(almost) Question~\ref{Q3}}} , e^{-st} \right\rangle \\
& = \underbrace{\ml{L}\left(\Cor_{[N(K_1)],\iota_V[N(K_2)]}\right)(s)}_{\text{Laplace transform}}  ,
\end{align*}
and 
\begin{align*}
\underbrace{\zeta_{\text{Conv}}(K_1,K_2,s)}_{\text{Question~\ref{Q1}}}& =\langle \underbrace{\mathsf{T}_{K_1,K_2}}_{\text{Question~\ref{Q2}}} , t^{-s} \rangle = \left\langle \underbrace{\Cor_{[N(K_1)],\iota_V[N(K_2)]}(t)}_{\text{(almost) Question~\ref{Q3}}} , t^{-s} \right\rangle \\
& = \underbrace{\M\left(\Cor_{[N(K_1)],\iota_V[N(K_2)]}\right)(s)}_{\text{Mellin transform}} .
\end{align*}

Unfortunately, the De Rham currents $[N(K_1)]$ and $\iota_V[N(K_2)]$ do not belong to the anisotropic (De Rham current version of the) spaces $\mathcal{H}^{M,N}(S\IT^d)$ defined in Section~\ref{s:sobobo}. Yet, up to application of a map of the form $(x,\theta)\mapsto (x+x_K(\theta),\theta)$ (where $x_K$ denotes the inverse of the Gauss map of $K$), they belong to these spaces and additional technical work is required to handle this extra transformation. Also, the key formula~\eqref{e:current-zeta} requires a geometric ``uniform transversality''  assumption, which is satisfied as a consequence of the strict convexity of the sets $K_1,K_2$.  We do not enter into these discussions, and instead refer the possible interested reader to the article~\cite{DangLeautaudRiviere21}.

\section*{Acknowledgements} We warmly thank B.~Chantraine, N.B.~Dang, F.~Faure, Y.~Guedes-Bonthonneau, D.~Han-Kwan and J.~Viola for interesting discussions related to this work. 
ML and GR are partially supported by the Agence Nationale de la Recherche through the PRC grant ADYCT (ANR-20-CE40-0017). GR also acknowledges the support of the Institut Universitaire de France and of the PRC grant ODA from the Agence Nationale de la Recherche (ANR-18-CE40-0020).

\bibliographystyle{alpha}
\bibliography{allbiblio}
\end{document}